\documentclass[11pt]{article}
\usepackage{geometry}                
\usepackage{graphicx}
\usepackage{amssymb}
\usepackage{epstopdf}
\usepackage{psfrag}
\usepackage{amsthm}
\def\Bbb{\mathbb}
\def\R{{\Bbb R}}
\def\C{{\Bbb C}}

\def\H{{\Bbb H}}

\def\mbf0{\mathbf{0}}

\newcommand{\bmat}{\begin{pmatrix}}
\newcommand{\emat}{\end{pmatrix}}

\theoremstyle{theorem}
\newtheorem{thm}{Theorem}
\theoremstyle{theorem}
\newtheorem{cor}{Corollary}
\theoremstyle{theorem}
\newtheorem{lem}{Lemma}
\theoremstyle{definition}
\newtheorem*{defn}{Definition}
\theoremstyle{remark}

\newtheorem*{rem}{Remark}

\title{ 
All finitely presentable groups from  \\
link complements
and Kleinian groups}

 \author{Iain R. Aitchison\\
 \\
\small{Dedicated to Professor Shinji Fukuhara on his 65th birthday} 
 }
 \date{}     

\begin{document}

 \maketitle

\begin{abstract}

Klein defined geometry in terms of invariance under groups actions; here we give a discrete (partial) converse of this, interpreting all (finitely presentable) groups in terms of the geometry of hyperbolic  3-manifolds (whose fundamental groups are, appropriately, Kleinian groups).  
For $G^*$ a Kleinian group of isometries of hyperbolic 3-space $\H^3$, with $M_{G^*}\cong \H^3/G^*$ a non-compact $N$-cusped
orientable   3-manifold of finite volume, let ${\mathcal P}_{G^*} \subset S^2_\infty=\partial\bar \H^3$ be its dense set of parabolic fixed points.  Let 
  $\bar M_{G^*} := \H^3 \cup {\mathcal P}_{G^*}/G^*$ be the   3-complex obtained by compactifying each cusp of $M_{G^*}$ with an additional point. This is the 3-dimensional analogue of the
  standard compactifcation of cusps of hyperbolic Riemann surfaces.
  We prove that every finitely presentable group $G$ is of the form
  $G = \pi_1(\bar M_{G^*})$, in infinitely many ways:
 thus every finitely presentable group arises as the fundamental group of an orientable  3-complex $\bar M$ -- 
 denoted as a `link-singular' 3-manifold --
 obtained from a hyperbolic link complement
by coning each boundary torus 
of the link  exterior 
    to a distinct point.

We define the  \emph{closed-link-genus},
$clg(G)$,  of any finitely presentable group $G$, which completely characterizes fundamental groups of closed orientable 3-manifolds:  $clg(G)=0$  if and only if $G$ is the fundamental group of a closed orientable 3-manifold. Moreover   $clg(G)$  gives an upper bound for the concept \emph{genus}($G$) of genus defined earlier  by Aitchison and Reeves,  and   in turn is bounded by the minimal number of 
relations among all finite presentations of $G$.  

Our results place some aspects of the study of finitely presentable groups more centrally  within both classical and modern 3-manifold topology: accordingly, proofs given are expressed in these terms, although
some can be seen naturally in 4-manifold topology.
\end{abstract}

 
 \section{Introduction and preliminaries}

 For general background, we refer the reader to Lyndon    and Schupp \cite{LS}, Hempel \cite{He}, Rolfsen \cite{Rol} and Thurston \cite{Th}.

 Presentations for groups, and their relationship with fundamental groups of topological complexes,
 has motivated much of the 20th century research into the classification of manifolds. The successful classification of the fundamental groups of surfaces naturally led to the desire to classify 3-manifold fundamental groups, and understand the extent to which uniqueness holds. In the case of simply-connected compact 3-manifolds, uniqueness is essentially equivalent to the Poincar\'e Conjecture, recently resolved by Perelman in his solution of Thurston's Geometrization Conjecture for 3-manifolds. In principle, fundamental groups of 
 compact 3-manifolds are now algorithmically classifiable (Bridson \cite{B1}).  It has been known for decades that every finitely presentable group does arise (non-uniquely)  as the fundamental group of some closed orientable 4-manifold.
Not all finitely presentable groups arise as the fundamental groups of compact 3-manifolds -- see for example,  Kawauchi \cite{Ka} or Shalen \cite{Sh} -- and moreover  the isomorphism problem for finitely presentable groups is algorithmically unsolvable (\cite{Ad, Rab}):
  it is of interest to make more precise the nature of the distinction between compact 3-manifold-groups and arbitrary finitely presentable groups.

   Quinn has shown that all finitely presentable groups arise
 as the  fundamental groups of  \emph{non-orientable} 3-complexes,  \emph{allowing boundary}, but where
all  vertex links are either spheres or  \emph{projective planes}. This construction is elucidated in \cite{H-A, H-AM}.   In this paper we show they arise as the  fundamental groups of  \emph{orientable} 3-complexes, with \emph{empty} boundary, with  
all  vertex links being either spheres or  \emph{tori}. 
As a consequence,  the combinatorial group theory of finitely presentable groups can be set more centrally in the theory of orientable 3-manifolds, geometric structures and the theory of knots and links.

    Every closed orientable surface can be triangulated;
after deleting vertices, taking the universal cover topologically produces the  ideal triangulation of the hyperbolic plane. Geometrically, we see every closed surface fundamental group arising 
from the compactification of a hyperbolic Riemann surface, with finitely-generated free fundamental group, by replacing cusp points.
  All closed (ie, compact with empty boundary) orientable topological 3-manifolds can be given an essentially unique piecewise-linear (`PL') structure -- and hence can be viewed as a union of tetrahedra with all faces identified in pairs. Every PL-3-manifold admits an essentially unique smoothing as a differentiable manifold. Similar arguments show that every compact orientable 3-complex, obtained by arbitrary pairwise identification of faces of a disjoint union of tetrahedra, admits an essentially unique differentiable structure in the complement of its vertices.
The following is well known (see for example \cite{He, Th}):

\begin{lem}
For 
a compact orientable triangulated 3-complex $M$ without boundary, obtained by pairwise face identifications of a finite number of tetrahedra, the following are equivalent:
\begin{enumerate}
\item The Euler characteristic $\chi (M)$ of $M$ vanishes: $\chi (M) = 0$;
\item $M$ is a 3-manifold;

 \item all interior vertex links are 2-spheres.
   \end{enumerate}

 \end{lem}

A  simple weakening of the concept of (compact, orientable) 3-manifold is to allow finitely many vertex links which are the next simplest closed orientable surface, that is, a torus. (In dimension 2, a link can only be  a circle.)
 Suppose $M$ is such a 3-complex, and we delete an open cone neighbourhood of each singular point. The result is  a 3-manifold $M^*$ with non-empty boundary a finite union of tori, and every such manifold can be obtained from a closed orientable 3-manifold $\bar M$ 
 by deleting an open tubular neighourhood of an embedded link ${\mathcal L}\subset {\bar M}$.
    Since coning a 2-sphere boundary component to a point does not change the fundamental group,  in the following we 
    generally consider compact orientable 3-manifolds    with no 2-spheres in their  
 boundary.

  \begin{rem}
  The Euler characteristic of a surface determines its possible constant-curvature geometries. Accordingly, spheres and tori 
 play a major role in the structure of 3-manifolds: embedded spheres arise from $\pi_2$, giving the prime decomposition of 3-manifolds; 
 embedded tori yield the JSJ decomposition of 3-manifolds (also mirrored in group-theoretic constructions); and any aspherical, atorioidal closed orientable 3-manifold with infinite fundamental group 
 admits a complete metric of negative curvature, according to Perelman-Thurston. 

 \end{rem}

 \begin{lem}
 Suppose $M$ is an arbitrary compact, orientable 3-manifold with 
 $\partial M\not= \emptyset$, but
 containing no 2-sphere components.
 The following are equivalent:
 \begin{enumerate}
\item $\partial M $ is a disjoint union of finitely many tori;
\item the Euler characteristic vanishes:  $\chi (M) = 0$. 
\end{enumerate}
\end{lem}

Let ${\mathcal M}^3_\chi $ denote the set of compact, connected, orientable 3-manifolds $M$ with   Euler characteristic  $\chi (M)=\chi$.
  Every finite-volume hyperbolic 3-manifold is uniquely the interior of some $M  \in {\mathcal M}^3_0 $.

 \begin{defn}
A {\it link-singular} 3-manifold is any compact orientable 3-complex obtained from a compact orientable 
3-manifold $M\in  {\mathcal M}^3_0$  by attaching cones to some or all  boundary tori:  
 all   interior vertex links are thus either spheres or tori.  
\end{defn}

For $M  \in {\mathcal M}^3_0 $, denote by $M_{C\partial}$   the orientable 3-complex obtained 
by coning each boundary torus to a distinct point (setting $M_{C\partial}\equiv M$ when $\partial M = \emptyset$).     
When ${\mathcal L} \subset M$ is a finite-component link in some  closed orientable 3-manifold $M$,
there is an open tubular neighbourhood $N{\mathcal L}$ giving the link exterior 
  $M_{\mathcal L}    :=  M -N{\mathcal L} \in {\mathcal M}^3_0$. 
 We let $  M_{C\mathcal L} := (M_{\mathcal L})_{C\partial}$ denote the link-singular 3-manifold obtained from
 $  M_{\mathcal L}$ by attaching a cone to each boundary torus of the exterior.

   \begin{rem}
Instead of coning boundary tori to points, we can attach a `solid torus' or `donut' -- $S^1\times D^2$ -- 
to some or all boundary tori, to obtain a closed orientable 3-manifold. 
 Any closed orientable 3-manifold $M$ can be obtained from any given one $M_0$  by Dehn surgery on a link ${\mathcal L}\subset M_0$, which essentially means the deletion of solid tori neighbourhoods of all link components, and reattaching these by a possibly homotopically non-trivial homeomorphism of their boundary tori, specified by an integer assigned to each link component. This follows from work 50 years ago by Rohlin, Lickorish and Wallace: Craggs and Kirby independently described the generators of the equivalence relation on the surgery instructions required to produce the same 3-manifold, notably in the case $M_0 = S^3$.
   \end{rem}

 We give  elementary knot-theoretic proofs of the following:

\begin{thm}\label{MT}
 Let $G$ denote an arbitrary finitely presentable group. Then
 $G\cong \pi_1(M_{C{\mathcal L}})$ for some link ${\mathcal L}$ in  a closed 3-manifold $M$, whose complement  
 $M- {\mathcal L}$ admits a complete  metric of constant curvature $-1$ and  finite volume. 
 Thus every finitely presentable group is the fundamental group of a link-singular 3-manifold, where the  complement of all singular points   admits a complete hyperbolic metric of finite volume.
 \end{thm}

    As a corollary, we obtain our main theorem, using the
 3-dimensional analogue of the
  standard compacitifcation of cusps of hyperbolic Riemann surfaces:
 For $G^*$ any Kleinian group of isometries of hyperbolic 3-space $\H^3$, with $M_{G^*}\cong \H^3/G^*$ a non-compact $N$-cusped
orientable   3-manifold of finite volume, let ${\mathcal P}_{G^*} \subset S^2_\infty=\partial\bar \H^3$ be its dense set of parabolic fixed points.  Let 
  $\bar M_{G^*} := \H^3 \cup {\mathcal P}_{G^*}/G^*$ be the   3-complex obtained by compactifying each cusp of $M_{G^*}$ with an additional point.

   \begin{thm} Every finitely presentable group $G$ is of the form
  $G = \pi_1(\bar M_{G^*})$, for  infinitely many   Kleinian groups $G^*$.
    \end{thm}

   \begin{rem}  If $G$ is the fundamental group of a closed orientable 3-manifold $M$, then $G$ can be realized as the fundamental group of a closed link-singular 3-manifold with any even number of singular points by $\#$-summing with some number of copies of $\Sigma T$, the suspension of a torus,  by the Seifert-van Kampen Theorem, since $\pi_1(\Sigma T) = 1$.  
    \end{rem}

   \begin{rem}
      That our main result  might be true was motivated by the 
 well-known result of Lickorish and Wallace that every  {closed} orientable 3-manifold can be obtained by surgery on some link in the 3-sphere (see \cite{Th}),
 that the connect sum $M_n :=  \#_nS^1\times S^2$ has fundamental group $\pi_1(M_n) \cong F_n$ (the   free group of rank $n$), and that all 
 finitely presentable groups arise as quotients of a free group. 
   \end{rem}

We prove the theorem by constructing a very simple link ${\mathcal L}_{P_G}$ in the 3-sphere $S^3$ from a given finite  presentation $P_G$  of a group $G$, analyzing the Wirtinger presentation of the fundamental group of its complement, and showing that killing off certain elements of the fundamental groups of peripheral tori yields the desired 3-complex with previously specified finite presentation of any arbitrary group. Similar links have arisen in related    applications in the past, and accordingly we introduce a convenient more symmetric refinement of these constructions (\cite{AkMa, Hu, Lic2, Lic3}).

 \section{Group presentations and link projections}

When referring to a finite presentation we will 
consider both the generating set and the set of defining relations as ordered sets.
Let 
 $P_G= \langle X_1,X_2,\dots , X_n \ | \ R_1, R_2,\dots,  R_k \rangle$
 be a finite (ordered) presentation for an arbitrary finitely presentable group $G$. Each relation 
 $R_j$ is thus a word $w_j(X_1,X_1^{-1}, \dots , X_n, X_n^{-1})$ in the monomials $X_i,\, X_i^{-1}$, and $G$ is obtained from the free group $F_n$ of rank $n$ with presentation 
 $P_n := \langle X_1,X_2,\dots , X_n \ | \  -  \rangle$
 by taking the quotient group $F_n/N_R$ where
 $N_R =  \langle\langle R_1,R_2,\dots , R_k   \rangle\rangle$
is the normal closure of the set of relations.

\begin{figure}[htbp]  
   \centering
   \includegraphics[width=4in]{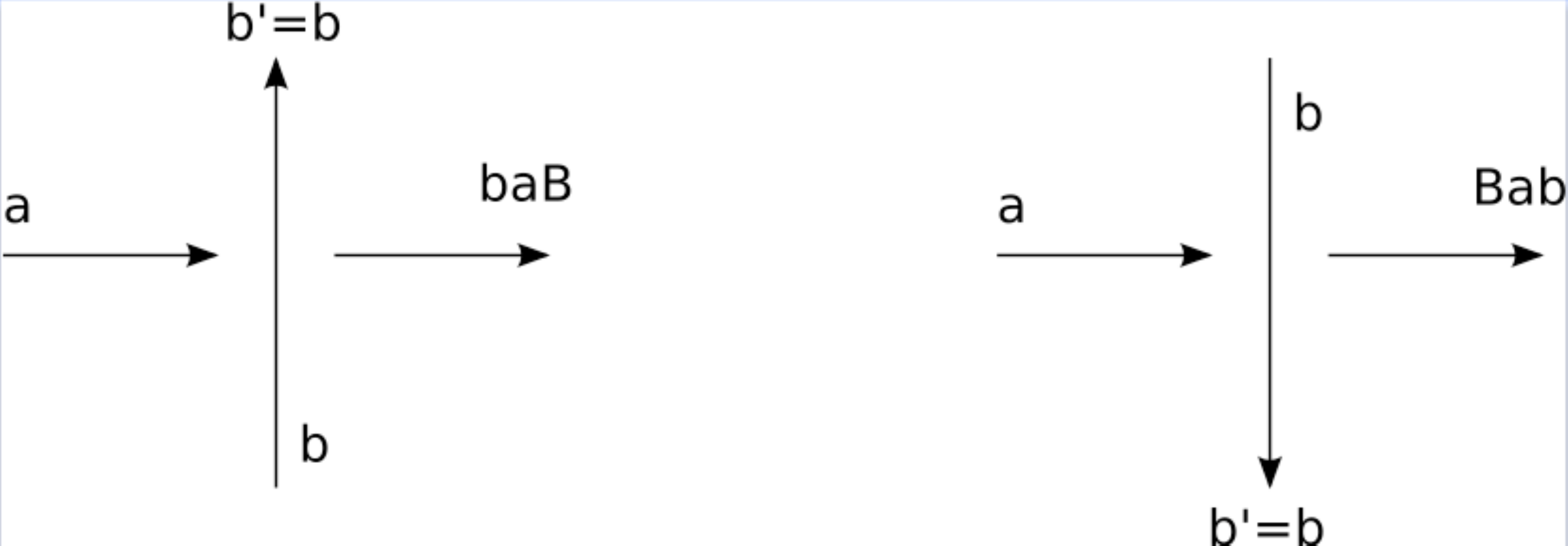} 
     \caption{ $a,b,b$' denote meridians of  arcs at a  `negative' crossing (left) and a  `positive'  crossing (right). Loops around arcs are oriented using the right-hand rule.  Capitals denote inverses: $B :=  b^{-1}$.}

   \label{fig:Wirtinger}
\end{figure}

We associate to $P_G$ a link in $S^3$ as follows: Each relation $R_j$ has a well defined length  $l_j = |w_j|$, and we set $L= \sum_1^k \, l_k$, the 
total length of the presentation. It is well known that the complement of a disjoint union of $t$ unknotted and unlinked circles has fundamental group $F_t$ free of rank $t$. From such a link, we will change crossings to introduce Wirtinger relations and produce an appropriate  Wirtinger presentation for the new link complement according to the conventions of Figure \ref{fig:Wirtinger}: labels correspond to generators of the fundamental group of the complement
assigned to each overarc of a link diagram.

 The standard Wirtinger presentation of the fundamental group of a link complement in the 3-sphere $S^3$ arising from a planar link projection with $N$ crossings 
 involves $N$ generators and $N$ relations. Generators correspond to small meridian loops linking each of the  $N$ oriented over-arcs of a projection, oriented according to the right-hand rule;
 the corresponding 4-valent planar graph has $N$ vertices and $2N$ arcs, which  combine to create the  $N$ overarcs.

 \bigskip

\begin{itemize}
\item Begin with $n+k$ concentric circles in the plane, labeled $g_1, \dots, g_n, r_1, \dots , r_k$ as indicated in Figure \ref{fig:concentriclabel2}, with anticlockwise orientation.
These are in correspondence with the $n$ generators and $k$ relations of $P_G$, with $g_i$ corresponding to $X_i$, $r_j$ to $R_j$. We refer to components of the link accordingly as \emph{generator components} and \emph{relation components}.

\begin{figure}[htbp] 
   \centering
   \includegraphics[width=2.5in]{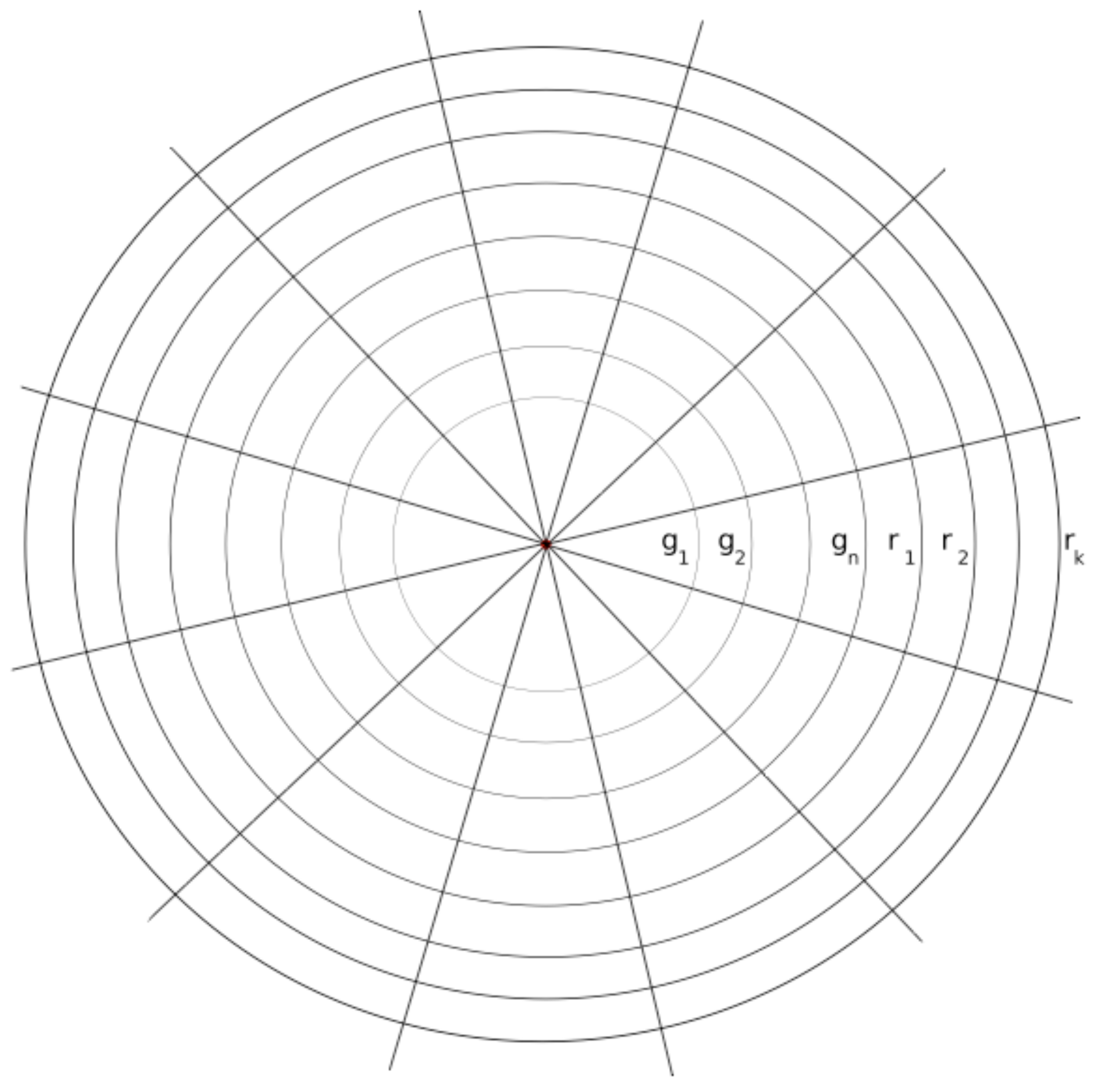}
  \caption{ 
  $n+k=8$ labelled concentric circles, each decomposed into $L=12$ arcs by rays from the origin.}
    \label{fig:concentriclabel2}
 \end{figure}

\item Concentrically divide each circle into $L$ equal length arcs. To each of the $L$ segments  there corresponds   a unique monomial
$ X_i^{\pm 1}$ occurring in a unique relation $R_j$ of $P_G$: as we read the monomials $X_i^{\pm 1}$  in the concatenated word
$R_1\cdots R_j\cdots R_k$ from left to right,  we  move anticlockwise from one segment to the next.

\item Each relation $R_j$ corresponds to $l_j$ consecutive segments of $r_j$. We will make $L$ modifications of the link, using the same procedure in each segment: if $X_i^\epsilon$ is associated to a given segment of $r_j$,
we create simple linking between the components arcs  of $g_i$ and $r_j$ within the segment, with linking number $\epsilon = \pm 1$.
This is illustrated in Figure \ref{fig:looplabelminus}, where we have conveniently substituted horizontal line segments for the concentric circular arc segments, labeled $g_1, \dots, g_n, r_1,\dots , r_k$ read top to bottom. Note that we maintain the labels assigned to components of the original unlink, and here add labels to overarcs of link components corresponding to their respective meridians in the fundamental group of the link complement.

\begin{figure}[htbp]  
   \centering
   \includegraphics[width=2.1in]{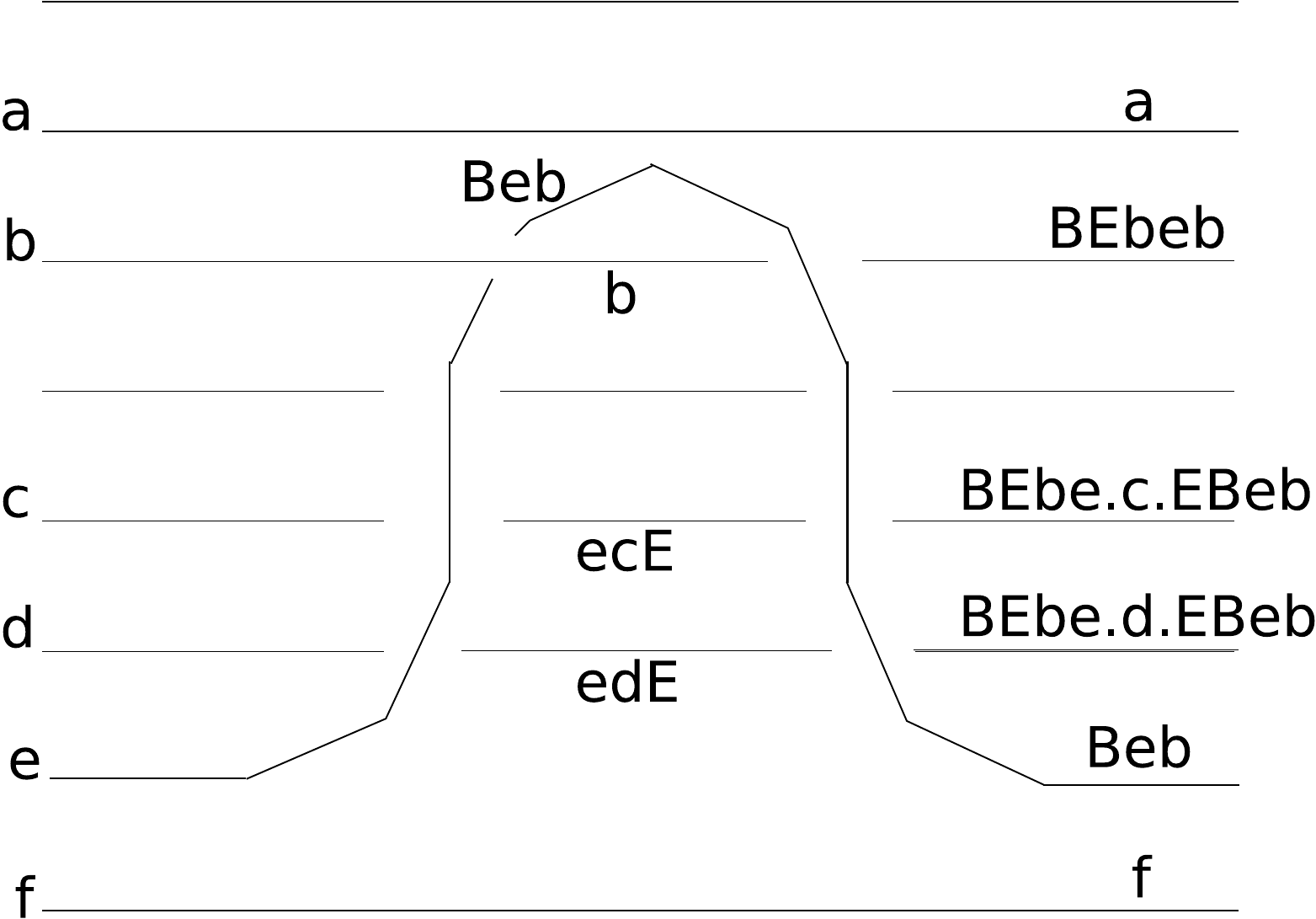} 
\qquad\qquad \qquad \qquad
   \includegraphics[width=2.2in]{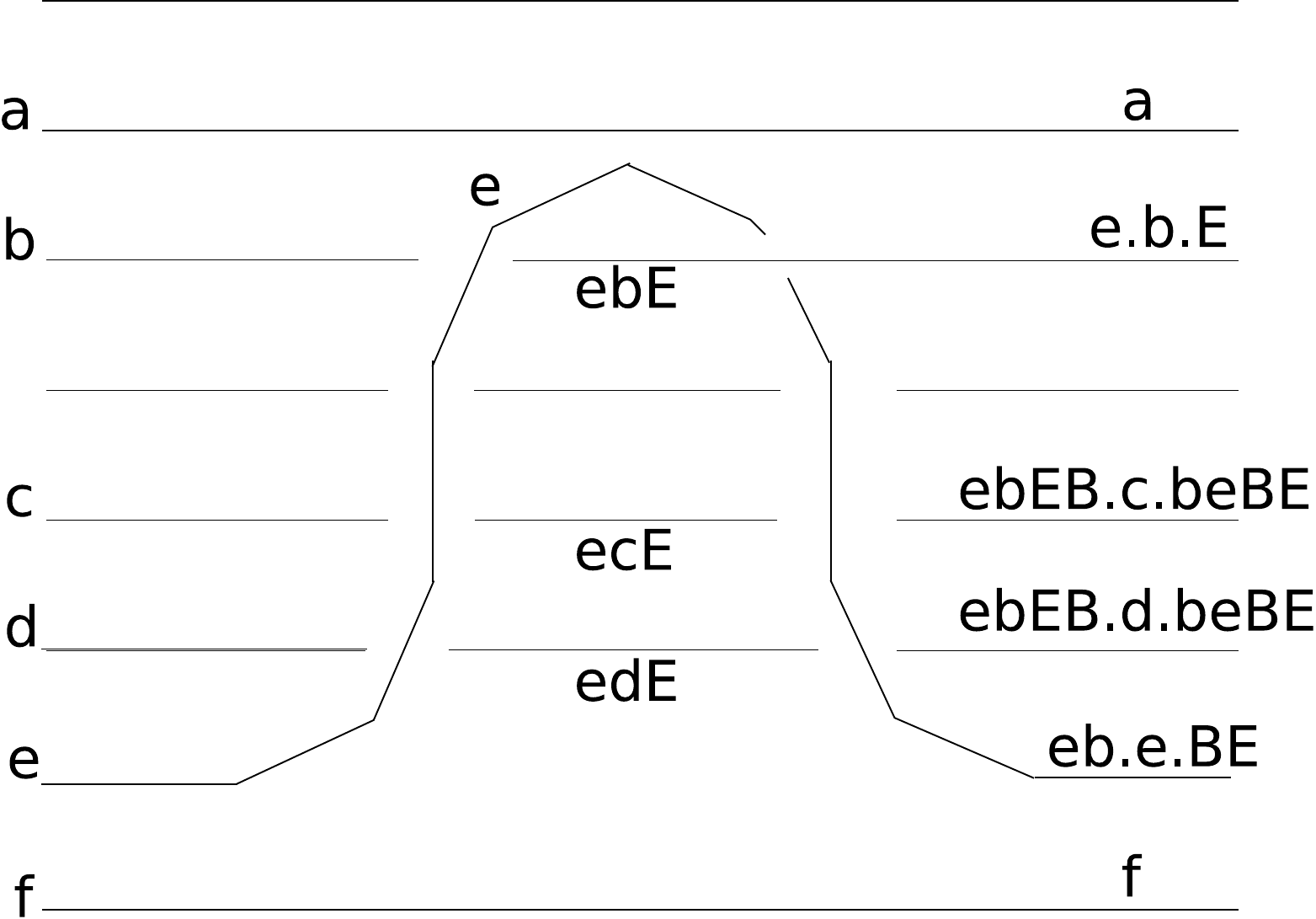} 
   \caption{Local linking from the monomials $X_i^{+1}$ (left) and $X_i^{- 1}$ (right)  in the relation $R_j$. The left arc of generator component $g_i$ has meridian label  $b\in \pi_1(S^3-{\mathcal L}_{P_G})$; the  relation component $r_j$ has left arc labeled $e\in \pi_1(S^3-{\mathcal L}_{P_G})$. Components are oriented left to right: $a,b,c$ correspond to meridians of generator component arcs; $d,e,f$ are meridians of relation component arcs.  Labels of arcs on the right of each figure arise from the convention for Wirtinger relations.}
   \label{fig:looplabelminus}
\end{figure}

\item  Having made $L$ modifications, we obtain a link 
${\mathcal L}_{P_G} = {\mathcal G}_{P_G}  \cup  {\mathcal R}_{P_G} $ 
of $n+k$ components, and $N$ crossings, which is naturally the union of two tangled unlinks 
$ {\mathcal G}_{P_G}\  {\mathcal R}_{P_G}$, respectively the generator and relation components.
The link ${\mathcal L}_{P_G} $ has
the following properties:
 
\subitem (a) Components remain  planar, unknotted, and parametrized by
the standard angles of the unit circle, giving a closed $(n+k)$-stranded pure braid;
\subitem (b)   Finitely many crossings occur, with at most one for any given angle;

 \vfill
 \pagebreak
\subitem (c)   The set of  plane-projected crossings  decompose all component circles $g_i,r_j$  into 
subarcs with corresponding meridians labeled either $g_{i,s} $ or $r_{j,t}$, where the
  second subscripts are consecutive along each oriented circle, beginning with the first-labeled subarc corresponding to the segment at standard angle $0\equiv2\pi$.  The total number of subarcs is $2N$, and at each crossing, the two subarcs of the overarc have equal meridian generators of $\pi_1(S^3-{\mathcal L}_{P_G})$ assigned to them. The Wirtinger relations determine how consecutive meridians assigned to underarcs are conjugate to each other;

\subitem  (d)  Crossings are ordered by their order of occurrence in  the concatenation of relations, $R_1\cdots R_j \cdots R_k$.

\item Generator and relation components may cross and link each other;   two relation components may cross but remain unlinked; but generator components do not cross each other. Collectively both the set of generator components, and the set of relation components, continue to form  unlinks of respectively $n$ and $k$ components,  with complements in $S^3$ having free fundamental group respectively of ranks $n$ and $k$.  There are thus six types of crossings of interest
in the   link projection, according to the sign of the crossing, and which of a generator or relation component creates the overarc.  

 \item  It is convenient to use a mild variation of the standard Wirtinger presentation, using $2N$ generators for symmetry purposes, and $2N$ relations, $N$ of which being equivalent to the standard relations, and the additional $N$ asserting equalities for the additional generators.

{\small 
\begin{table}[htdp]
\caption{Wirtinger meridian data for six crossing types involving oriented components $g_i,\, r_j, \, r_m$}
\begin{center}
\begin{tabular}{|c||c|c|c|c|c|c|}
\hline
&&&&&\\

 $g_{i,s}$ over $r_{j,t}$ & $g_{i,s}$ over $r_{j,t}$ & $r_{j,t}$ over $g_{i,s}$ & $r_{j,t}$ over $g_{i,s}$& $r_{j,t}$ over $r_{m,u}$ &$r_{j,t}$ over $r_{m,u}$ \\
&&&&&\\
$+:\ A_{i,s;j,t}$&$- : \ B_{i,s;j,t}$&$ +:\ C_{j,t;i,s}$&$ -:\  D_{j,t;i,s}$&$ +:\ E_{j,t;m,u}$&$ -:\ F_{j,t;m,u}$\\
&&&&&\\
\hline
\hline
&&&&&\\
$ g_{i,s+1} =  g_{i,s}  $ & $ g_{i,s+1} =  g_{i,s}  $ & 
$ r_{j,t+1} =  r_{j,t}  $ & $ r_{j,t+1} =  r_{j,t}  $ & 
$ r_{j,t+1} =  r_{j,t}  $ & $ r_{j,t+1} =  r_{j,t}   $ \\
&&&&&\\
$r_{j, t+1} =  $ & 
$r_{j, t+1} =   $& 
$g_{i,s+1} =   $ & 
$g_{i,s+1} =  $&
$r_{m, u+1} = $ & 
$r_{m,u+1} =  $\\ 
$ g_{i,s}^{-1} r_{j, t}  g_{i,s} $ &$g_{i,s} r_{j, t}  g_{i,s}^{-1}   $&
$r_{j,t}^{-1} g_{i,s}  r_{j,t}   $&$ r_{j,t} g_{i,s}  r_{j,t}^{-1}   $&$  r_{j,t}^{-1} r_{m,u}  r_{j,t}  $&
$ r_{j,t} r_{m,u}  r_{j,t}^{-1}  $\\
&&&&&\\
\hline
\end{tabular}
\end{center}
\label{Wirtinger}
\end{table}%
}

\item For any given component of any oriented link in $S^3$, corresponding generators for the (generalized) Wirtinger presentation  are all conjugate to each other. 
Table \ref{Wirtinger} enables us to record how meridian labels for arcs for generator and relation components are explicitly conjugate to each other, and thus to record contributions to words expressing longitudes -- planar parallels for any given component --
as products of conjugates of meridians. By potential abuse of notation, the labelling of arcs by meridians $g_{i,j}, r_{j,t}$ reflects the naming of components $g_i,r_j$.

 \end{itemize}

 \bigskip
 
 Assembling the data from each crossing, we obtain:

 \begin{lem}
 The defining Wirtinger presentation of $\pi_1(S^3-{\mathcal L_{P_G}})$ can be expressed  
$$  \pi_1(S^3 - {\mathcal L_{P_G}}) : =
 \langle
  g_{i,s} ,\ 
 r_{j,t} \, | \, 
 A_{i,s;j,t},\,  B_{i,s;j,t},\, C_{j,t;i,s},\, D_{j,t;i,s}, \, E_{j,t;m,u},\, F_{j,t;m,u} \\ \rangle
 $$
 where each symbol $X_{a,b;c,d}$ denotes the appropriate pair of relations for a crossing read off from   Table \ref{Wirtinger}. 
 
\end{lem}

 \section{Proof of Theorem \ref{MT}}
 
 To complete the proof, we identify longitudes for each component, and a distinguished meridian, corresponding to a choice of generators for the  
 fundamental groups of peripheral tori.
 The use of Tietze transformations allows us to 
 identify the quotient groups obtained by adding relations trivializing appropriate  meridians and longitudes, to obtain the initially given presentation $P_G$.
Accordingly we add generators  to the presentation corresponding to a choice of meridian and longitude for each peripheral torus, and relations expressing these in terms of existing meridian generators for the link complement.

 \bigskip

 \noindent
 {\bf Meridian generators:}
To each subset of generators $\{g_{i,s} \}_{s=1}^{s_i}$, there corresponds a new generator $X_i$ and defining relation 
$X_i  := g_{i,1}$; similarly, to $\{ r_{j,t} \}_{t=1}^{t_j}$ there corresponds a new generator $\bar R_j$ and defining relation $\bar R_j := r_{j,1}$. 
 The Wirtinger presentation of $\pi_1(S^3-{\mathcal L_{P_G}})$ can be modified by Tietze transformations to 
$  \pi_1(S^3 - {\mathcal L_{P_G}})  :=$
$$ \langle
  g_{i,s} ,\ 
 r_{j,t}, \, X_i ,\, \bar R_j | \, 
 A_{i,s;j,t},\,  B_{i,s;j,t},\, C_{j,t;i,s},\, D_{j,t;i,s}, \, E_{j,t;m,u},\, F_{j,t;m,u} , \, X_ig_{i,1}^{-1}, \,   \bar R_jr_{j,1}^{-1}\\ \rangle
 $$

 \bigskip
 
  \noindent
 {\bf 
Longitudes generators:} Each component is planar, unknotted, and hence admits a planar-parallel longitude in $S^3 -  {\mathcal L_{P_G}}$. Components are of two kinds --  corresponding to 
generators and relations of the given  presentation -- which are to be respectively surgered using the framing of a longitude, or coned so that meridians and longitudes become trivial.
The natural construction of the link leads to a simple product description for longitudes, corresponding to the product structure of concatenated relations given by  the presentation. A longitude for a component $g_i$ will be denoted $\gamma_i$, and for 
a component $r_j$, we denote a longitude by $\rho_j$.

To express $\gamma_i,\ \rho_j$ as  words in $\{ g_{i,s},\ r_{j,t}\}$, we follow the longitudes around the link, consecutively recording each 
oriented meridian generator corresponding to overarcs where the longitude passes under a component of ${\mathcal L}_{P_G}$.  Each longitude is thus a product of $L$ subwords,  corresponding to
the occurrences of   original generators $X_i^{\pm 1}$ in the  relation-concatenation $R_1\cdots R_j\cdots R_k$.    
Using the Wirtinger relations, as in Figure \ref{fig:looplabelminus}, we see that the involved arc of the longitudes parallel to components labeled $a,b,c,d,e,f$ contributes subwords
given by Table  \ref{long+} and Table \ref{long-}, with capitals again denoting inverses.

\begin{table}[htdp]
\caption{Wirtinger longitude subword contributions  from the occurrence of $X_i$ in $ R_j$}
\begin{center}
\begin{tabular}{|c||c|c|c|c|c|c|}
\hline
&\qquad\qquad&\qquad\qquad&\qquad\qquad&\qquad\qquad&\qquad\qquad&\qquad\qquad\\
{\bf  Component} & a   &b           &c       &d &         e&    f   \\
&&&&&&\\
\hline
&&&&&&\\
{\bf  Contribution} & 1  &Beb  &E.Beb   & E.Beb & b  &     1     \\
&&&&&&\\
{\bf  If $d,e,f=1$} & 1  &1    & 1 & 1 & $b=X_i$ &     1     \\
&&&&&&\\
\hline
\end{tabular}
\end{center}
\label{long+}
\end{table}

\begin{table}[htdp]
\caption{Wirtinger longitude subword contributions  from the occurrence of $X_i^{-1}
$ in $R_j$}
\begin{center}
\begin{tabular}{|c||c|c|c|c|c|c|}
\hline
&\qquad\qquad&\qquad\qquad&\qquad\qquad&\qquad\qquad&\qquad\qquad&\qquad\qquad\\
{\bf  Component} & a   &b           &c       &d &         e&    f   \\
&&&&&&\\
\hline
&&&&&&\\
{\bf  Contribution} & 1  &E    & E.ebeBE & E.ebeBE & EBe &     1     \\
&&&&&&\\
{\bf  If $d,e,f=1$} & 1  &1     & 1 & 1 & $B=X_i^{-1}$ &     1     \\
&&&&&&\\
\hline
\end{tabular}
\end{center}
\label{long-}
\end{table}%

 \bigskip
 
 \bigskip

Let $  M_{{\mathcal L}_{P_G}}$ denote the compact exterior of the disjoint union of open solid donut neighbourhoods of all link components
of ${\mathcal L}_{P_G}$, and let $  M_{C{\mathcal R}_{P_G}}$ denote the   link-singular 3-manifold obtained by
 attaching a cone to each boundary torus of relation components. Thus all meridian generators $r_{j,t}$ and longitudes $\rho_j$ of components
 $r_j$ are set equal to 1 as relations  added to the presentation of $  \pi_1(S^3 - {\mathcal L_{P_G}}) $, by the Seifert-van Kampen Theorem.

 \begin{lem}The group $\pi_1(  M_{C{\mathcal R}_{P_G}})$ admits a finite presentation $ :=
 \langle
   \, X_i\, | \, 
       \rho_j \rangle
 $.
 \end{lem}
 
 \noindent
 {\it Proof:} We apply Tietze transformations to the presentation of $  \pi_1(  M_{C{\mathcal R}_{P_G}}) :=$
 $$
 \langle
  g_{i,s} ,\ 
 r_{j,t}, \, X_i ,\, \bar R_j | \, 
 A_{i,s;j,t},\,  B_{i,s;j,t},\, C_{j,t;i,s},\, D_{j,t;i,s}, \, E_{j,t;m,u},\, F_{j,t;m,u} , \, X_ig_{i,1}^{-1}, \, , \, \bar R_j r_{j,1}^{-1}, \, r_{j,t},\, \rho_j \rangle
 $$
 $$ :=
 \langle
  g_{i,s} ,
 \, X_i ,\, \bar R_j | \, 
 g_{i,s+1}=g_{i,s},\,    X_ig_{i,1}^{-1},  \, \bar R_j ,\, \rho_j \rangle
   :=
 \langle
 X_i ,\, \bar R_j | \, 
    \bar R_j,\, \rho_j\ \rangle
     :=
 \langle
 X_i  \, | \, 
    \rho_j\\ \rangle. 
 $$ 
 
   We now identify the words defining each $\rho_j$ in $\pi_1(  M_{C{\mathcal R}_{P_G}})$, using the tables. From these,  we have:
 
 \begin{lem}
 (a) The group $\pi_1(  M_{C{\mathcal R}_{P_G}})$ admits a finite presentation $ :=
 \langle
   \, X_i\, | \, 
      R_j\rangle
 $.  (b) In  $\pi_1(   M_{C{\mathcal R}_{P_G}})$, $\gamma_i := 1$.
 
 \end{lem}

 Now perform 0-surgery on each generator component of ${\mathcal L}_{P_G}$, to obtain a closed link-singular 3-manifold $M_{P_G}$. This is achieved by attaching a solid donut $S^1\times D^2_i$
to each of the remaining boundary tori of $  M_{C{\mathcal R}_{P_G}}$, with the boundary of a  $D_i$  attached with framing determined by  
$\gamma_i = 1\in \pi_1(   M_{C{\mathcal R}_{P_G}})$. These elements are already trivial, and we obtain:
 
 \begin{thm}
 The fundamental group $ \pi_1(  M_{P_G})$ has presentation
 $ \pi_1(  M_{P_G}) :=  \langle
   \, X_i\, | \, 
      R_j\rangle.$ 
 
 \end{thm}
 
   \begin{cor}
   Every finitely presentable  group admits a representation as the fundamental group of a closed, orientable link-singular 3-manifold,
 obtained by adding cones to the toral boundary components of a link exterior in some closed orientable 3-manifold. 
    \end{cor}
    
   \noindent
 {\it Proof:} Reverse the order of attaching cones and solid donuts to $S^3 - {\mathcal L}_{P_G}$. Cones are then attached to boundary tori of a link exterior in the connect-sum $\#_n\, S^1\times S^2$, the result of 0-surgery on the unlink ${\mathcal G}_{P_G}$ of $n$ components.

 \section{Genus: characterising compact  orientable 3-manifold groups}

 In \cite{AR}, it is shown by a different construction that every finitely presentable group arises as the fundamental group of a `singular 3-manifold', in infinitely many ways: this is constructed  from an orientable compact 3-manifold by coning one boundary component to a point, producing the only non-manifold (singular) point.
 Note that the resulting singular 3-manifold   may have additional boundary components.     
 
   \begin{defn} (Aitchison and Reeves \cite{AR})
 The \emph{genus} of $G$, denoted \emph{genus}$(G)$,  is defined as the least possible genus of such a boundary component for any such singular 3-manifold realizing the group $G$.  
   \end{defn}

   \begin{rem}
 A question raised in \cite{AR} 
 was whether or not every group could be realized without additional boundary components. We have answered this in the affirmative
 above.
    \end{rem}

   \begin{defn}
    The \emph{closed link  genus}, denoted $cl g(G)$,  of a finitely presentable group $G$ is the 
 mininimal number of boundary tori among all $M \in {\mathcal M}^3_0$ such that $G\cong \pi_1(M_{C\partial})$.
    \end{defn}

   \begin{rem}
   By  Theorem \ref{MT}, this is well defined for any finitely presentable group $G$.  
    \end{rem}

\begin{thm}
A finitely presentable group $G$  has $clg(G)= 0$ if and only if $G$ is the fundamental group of a closed orientable 3-manifold.
\end{thm}

\noindent
{\it Proof:}  If $G\cong \pi_1(M)$ for some closed orientable 3-manifold $M$, then $clg(G) = 0$ follows from the definition.
If $clg(G) = 0$, then $G$ is realized by a closed link-singular 3-manifold $M$  with no singular points, which is thus a closed 
orientable 3-manifold.

\bigskip

If $G$ is the fundamental group of a compact  orientable 3-manifold $M$ with a   boundary component of positive genus, then $G$ is also the fundamental group of some  closed link-singular 3-manifold $M^*_{C\partial},\ M^*\in {\mathcal M}^3_0$, by the main theorem, generally not 
 free of singularities (since otherwise $G$ is also a closed 3-manifold group). Conversely, by allowing additional boundary components, 
 it is conceivable that we may decrease the number of singular points required to realize $G$. 
  
\begin{thm}
Suppose $G$ admits a finite presentation with $k$ relations. Then
\ $$
\emph{genus}(G) \leq clg(G) \leq k.
$$
\end{thm}

\noindent
{\it Proof:} By the construction above, we know there exists a closed link-singular 3-manifold with at most $k$ 
singular points which realizes $G$. Thus $clg(G)\leq k$. 
In \cite{AR}, it is shown  that if there exists a singular 3-manifold $M$ realizing $G$, with two singular points respectively arising from cones on surfaces of genus $g_1,g_2$, then there is a singular 3-manifold $M'$ realizing $G$ with a single cone on a surface 
of genus $g_1+g_2$ replacing these. Thus if there are $k$ cones on tori in a closed link-singular 3-manifold $M$ realizing $G$,
there is a singular 3-manifold $M'$ realizing $G$ with exactly 1 singular point, arising from a surface of genus $k$. Thus 
$\emph{genus}(G)\leq clg(G)$.

 \bigskip
 
 \noindent
 {\bf Example.} Baumslag--Solitar groups with presentations of form 
 $G:=\langle x_1, x_2  \, | \, x_1x_2^mx_1^{-1}x_2^{-n}\rangle$ were shown in \cite{AR} to have 
 \emph{genus}$(G)= 1$ when $|n|,\ |m|\not=1$. Since they are 1-relator groups, we also have $clg(G) = 1$. It is well known that such groups are not fundamental groups of compact orientable 3-manifolds,
 with or without boundary \cite{Sh}.  

\section{Hyperbolic link exteriors suffice}

Every non-compact hyperbolic 3-manifold $M$ of finite volume is the interior of 
a compact 3-manifold $M\in {\mathcal M}^3_0$  with nonempty boundary a disjoint union of tori.  We  show that every finitely presentable group $G$ can be realized by coning boundary components of the exterior of a link whose complement admitting a hyperbolic structure.
  
\begin{thm}
Suppose $G$ is an arbitrary finitely presentable group. Then there are infinitely many 
$M\in {\mathcal M}^3_0$ with interior admitting a (unique)  hyperbolic  structure, such that $G\cong \pi_1(M_{C\partial})$.
\end{thm}

\noindent
{\it Proof:} 
We recall the following theorem of Myers, quoted from \cite{My1}:

\begin{thm}
{\rm (Myers)} {Let $M$ be a compact, connected 3-manifold. Suppose $J$ is a compact (but
not necessarily connected), properly embedded 1-manifold in $M$. $J$ is homotopic rel $\partial J$
to an excellent 1-manifold $K$ if and only if $J$ meets every 2-sphere in $\partial M$ in at least
two points and every projective plane in $\partial M$ in at least one point. In this case there
are infinitely many such $K$ with nonhomeomorphic exteriors. Moreover, each $K$ can be
chosen so that it is ribbon concordant to $J$}.

\end{thm}

An excellent 1-manifold is one whose complement is an excellent 3-manifold:
Excellent 3-manifolds admit
hyperbolic structures, i.e., Riemannian metrics on their interiors having constant
sectional curvature -1. Heuristically, homotoping $J$ to $K$ punctures all homotopically essential spheres and tori without creating new ones (other than peripheral tori) in its complement.  

 \medskip
 
To conclude the proof, consider 
${\mathcal  L }_{P_G}  = {\mathcal  G}_{P_G} \cup {\mathcal  R }_{P_G}$.
Carry out 0-surgery on the sublink of generating  components  
${\mathcal  G }_{P_G} $, so that ${\mathcal  R }_{P_G}$ becomes a new link $ J := {\mathcal  R' }$  of relation components in $M =\#_n\, S^2\times S^1$.   
Let $K = {\mathcal  R'' }$ be any excellent 1-manifold obtained from $J$ by Myers' Theorem. Since $ {\mathcal  R'' }$ is homotopic to
$ {\mathcal  R' }$, after coning we find $\pi_1(M_{C{\mathcal  R' }}) \cong \pi_1(M_{C{\mathcal  R'' }})$: this follows since meridians of relation components are set to $1$, and so Wirtinger generators 
and relations corresponding to  crossings of relation components are irrelevant to the presentation of these groups. We may change such crossings at will, corresponding to homotopy, without changing the resulting groups. 
 
  \bigskip
  
  For $G^*$ a Kleinian group of isometries of hyperbolic 3-space $\H^3$, with $M_{G^*}\cong \H^3/G^*$ a non-compact $N$-cusped
orientable   3-manifold of finite volume, let ${\mathcal P}_{G^*} \subset S^2_\infty=\partial\bar \H^3$ be its dense set of parabolic fixed points.  Let 
  $\bar M_{G^*} := \H^3 \cup {\mathcal P}_{G^*}/G^*$ be the   3-complex obtained by compactifying each cusp of $M_{G^*}$ with an additional point. This is the 3-dimensional analogue of the
  standard compactifcation of cusps of hyperbolic Riemann surfaces. As  corollaries to the last theorem, we have:
  
  \begin{thm}
Every finitely presentable group $G$ is of the form
  $G = \pi_1(\bar M_{G^*})$, in infinitely many ways.
  \end{thm}
 
\begin{cor}
Any invariant of non-compact hyperbolic 3-manifolds of finite volume, whose values can be  ordered, defines an invariant of 
finitely presentable groups, by minimizing values over all hyperbolic link complements for which the addition of cones to their exterior realizes any given group   $G$.
\end{cor}
 As an example:
\begin{defn}
The \emph{volume} \emph{vol}($G$) of a finitely presentable group $G$ is the least volume of any hyperbolic link complement ${\mathcal L}\subset M$ such that $M_{C\mathcal L}$ realizes $G$.
 \end{defn}

\section{Concluding remarks}

\begin{enumerate}

\item \emph{Volume:} 
 The trivial group arises by coning the boundary tori of the exterior of any link with hyperbolic complement in $S^3$, since the fundamental group 
 of the complement is generated by meridians. When $G$ is the trivial group,  \emph{vol}($G$) is of course bounded by the smallest volume of any hyperbolic link complement in $S^3$. What can be said of the groups realized by coning on the smallest volume hyperbolic  link complements in $\#_n\, S^2\times S^1$?

\item \emph{Relative hyperbolicity:} Gromov's definition of relatively hyperbolic groups is motivated by the fundamental groups
and peripheral subgroups of hyperbolic link complements. The construction of the previous section shows the naturality of this concept 
and   the direct use of the hyperbolic geometry of link complements to understand finitely presentable groups, arising by Dehn surgery on link complements -- see for example Lackenby, and Fujiwara and Manning
\cite{La, FuMa}.

  \item \emph{Energy concepts:} A simple `generalization' of ${\mathcal L}_{P_G}$ is a tangle 
${\mathcal L} = {\mathcal G} \cup {\mathcal R} $
  of two unlinks 
  with respectively $n$ and $k$ components. Perform 0-surgery on components of ${\mathcal G } $, and cone on  components of  $ {\mathcal R} $. 
  The resulting  fundamental group can be read off from  the linking data of each component of $ {\mathcal R} $
  with components of $ {\mathcal G}$, and so this construction is   essentially  equivalent  to a homotopy of the components of $ {\mathcal R}_{P_G} $ in the complement of ${\mathcal G}_{P_G}  $.
  It would be interesting to use some notions of knot/link energies to find some canonical,  or a finite number of possible  links,
 representing a given group $G$  by minimizing some `relative' link-homotopy energy associated to the linking of the unlinks 
 $ {\mathcal G}, \, {\mathcal R}$, derived from such as the   original one first introduced by 
Fukuhara \cite{Fu1}. Given a presentation for a group $G$, there should be only finitely many possibilities, and generically perhaps one.
When is such a link hyperbolic, and in such cases, what is its volume compared to $vol(G)$?

\item \emph{Braids:} Our construction of ${\mathcal L}_{P_G}$
produces a link which is a pure braid: for example, the
Baumslag--Solitar groups with presentations of form 
 $G:=\langle x_1, x_2  \, | \, x_1x_2^mx_1^{-1}x_2^{-n}\rangle$
 give 3-braids which are of very simple form: 
 Recall that link is called \emph{quasipositive} \cite{Ru} if it is the closure of a braid which is the product of conjugates of the Artin generators $\sigma_i^{+1}$. This concept has its origins in links of singularities of algebraic curves in $\C^2\cong \R^4$, locally related to cones on 
 knots and links. It would be   interesting to relate these concepts with ${\mathcal L}_{P_G}$ for a given $P_G$, which 
 is the product of conjugates of {squares}
 $\sigma_i^{\pm 2}$, when exponents are positive: such is the case when $m,-n >1$ above. For $n,m >1$, the braids are alternating, and are thus fibered links with hyperbolic complements. Is it possible to find tangled unlink representatives for all finitely presentable groups, which are pure braids, so that their complements, before coning or surgery, are fibered and hyperbolic?
 In another direction, braid groups themselves are finitely presentable fundamental groups of hyperplane complements, and have been recently considered in the context of cryptography. Further comments on decision problems are given below.
 
\item \emph{4-dimensional considerations:}
Invariants of finitely presentable groups have naturally  been considered in the context of invariants of 4-manifolds: see for example Kotschick \cite{Ko3,Ko4}.
  As an example of an invariant for a finitely-presentable group $G$, the Hausmann-Weinberger invariant \cite{HaWe}  is defined as the minimal Euler characteristic $q(G)$ of a closed orientable 4-manifold $M$ with fundamental group $G$: for a given $M$, the Euler characteristic $\chi(M)$  is easily computed from a link description of $M$, but generally such calculations only give upper bounds for such invariants. Our invariants can be considered as a lower-dimensional `analogue' of this.
In this note we proved all results within the realm of classical 3-dimensional topology: a more  natural setting for some is in 
  the theory of 4-manifolds, which will be discussed in a subsequent note. 
  \item \emph{`Geography':} For a fixed finitely presentable group $G$,
  characterize the set of hyperbolic link complements 
 corresponding to $G$.

  \item \emph{Recent quantum and other invariants:} The results of this paper lead naturally to the reinterpretation of
  recent 3-manifold invariants in the context of combinatorial and geometric group theory, and in   `virtual' knots  and links. 
  
\item \emph{Number theory:} 
 We see that our construction makes contact with classical constructions from Riemann surface theory and the number theory of hyperbolic geometry. For hyperbolic Riemann surfaces, compactification of cusps is natural,
and relates number theory with parabolic fixed points  added to the circle at infinity of the hyperbolic plane: this construction is fundamental, pertaining  to `Monstrous Moonshine', and the Taniyama--Shimura--Weil Conjecture.
 Characterize finitely presentable groups arising from  \emph{arithmetic} link complements -- see Maclachlan and Reid \cite{MR}, and 
Neumann and Reid \cite{NR}.

\item \emph{Decision problems:} The following well known theorems raise the problem of identifying certain  hyperbolic link complements which are
 potentially very interesting.  Find `nice'  representative link complements for finitely presented groups having the properties of the following theorems:
For a survey of decision problems in group theory, see \cite{AD, LS}.
 \begin{thm} {\rm  (Higman \cite{Hig})}
 There is a finitely presented group $F$ that is universal for
  all finitely presented groups. This means that for any finitely presented 
 group $G$ there
is a subgroup of $F$ isomorphic to $G$.
  \end{thm} 
 There is an analogue -- \emph{ universal links} -- in 3-manifold topology:   the arithmetic Borromean link  ${\mathcal B}\subset S^3$ is  universal \cite{HLM}, as the set of  its branched-covers includes 
all closed orientable 3-manifolds. 
 Suppose $F\cong \pi_1(M^F_{C\partial}),\, M^F\in {\mathcal M}^3_0$, and $G$ is a subgroup of $F$. Then there is a covering $\pi : M^G_{C\partial} \to M^F_{C\partial}$,
 with $G\cong \pi_1(M^G_{C\partial})$, inducing a covering 
 $\pi : M^G := \pi^{-1}(M^F) \to M^F$. A torus boundary component of $M^F$ may be multiply covered, but restricting to any connected component of its preimage in $M^G$, the covering is a homeomorphism: such a 3-manifold $M^F$ has a very rich subclass of covering spaces of this form.

In the `opposite direction', we recall a  
 theorem   proved by C. F. Miller III:
 \begin{thm} {\rm (Miller \cite{Mi})}
There is a finitely presented group all of whose non-trivial quotient groups have insoluble word problem.
 \end{thm}  
Reconcile decision algorithms based on Tietze transformations of presentations with the recent proof of Thurston's Geometrization Conjecture: note 
trisection of an angle \emph{is} possible using ruler and compass constructions, \emph{allowing an additional marked point}: It is perhaps fruitful to now construe finitely presentable groups as fundamental groups
of link-singular manifolds, and define them as such rather then via presentations.

\item \emph{Presentation link calculus:} For surgery representations of 3-manifolds, stabilization and 
handle-sliding generate the equivalence relation between framed links
yielding the same 3-manifold. This creates an analogue of Tietze moves generating the equivalence relation between presentations of the same group: In the coned-link representation of finitely presentable groups, we can handle slide generating components over each other, and similarly relation components over each other. Relation components can slide over generating components, but not vice-versa. We can stabilize by 
adding  a  Hopf link separated from other components, consisting of a generator and a relation component, or add a new relation component,
unkotted and separated from the original link. Relation components can be homotoped arbitrarily in the complement of generating components, which is equivalent to homotopy in a connect-sum $\#_n\, S^2\times S^1$ obtained by 0-surgery on generator components: observe that (a)  the complement in $S^3$ of an $n$-component unlink has fundamental group which is free of rank $n$, as is 
$\pi_1(\#_n\, S^2\times S^1)$; and (b) sliding a relation component over a 0-surgered generator component yields a new link obtainable by homotopy of the relation component in the complement of the generator components, since all generator components are unknotted and unlinked.

\end{enumerate}

\section{Ackowledgements} The author would like to express thanks to Lawrence Reeves, Shinji Fukuhara, Makoto Sakuma, Akio Kawauchi, Brian Bowditch, Martin Bridson, Marc Lackenby, Dmitri Panov and 
Danny Calegari
for support, helpful conversations and encouragement, and Naoko and  Seichi  Kamada additionally for writing \cite{KaKa}, which inspired this work.

\bigskip

\noindent
Email: {iain at ms.unimelb.edu.au}
\end{document}